\documentclass[10pt,reqno]{amsart}

\usepackage{amscd,amssymb}
\usepackage[mathscr]{eucal}
\usepackage[all]{xy}
\usepackage{mathrsfs}
\usepackage{xypic}
\usepackage{amsfonts}
\usepackage{amsmath}
\usepackage{amsthm}
\usepackage{latexsym}
\usepackage{eufrak}

\theoremstyle{plain}
\newtheorem{theorem}[equation]{Theorem}
\newtheorem{corollary}[equation]{Corollary}
\newtheorem{proposition}[equation]{Proposition}
\newtheorem{lemma}[equation]{Lemma}

\theoremstyle{definition}

\theoremstyle{remark}

\newtheorem{remark}[equation]{Remark}

\newcommand\be[1]{\begin{equation}\label{#1}}
\newcommand\ee{\end{equation}}

\def\op{\operatorname}
\newcommand{\map}{\longrightarrow}

\def\k{{\Bbbk}}
\newcommand{\X}{\mathbb{X}}
\newcommand{\Y}{\mathbb{Y}}

\newcommand{\C}{\mathbb{C}}
\newcommand{\N}{\mathbb{N}}
\newcommand{\Z}{\mathbb{Z}}

\newcommand\lie[1]{\mathfrak{#1}}

\renewcommand{\b}{\lie{b}}
\newcommand{\g}{\lie{g}}

\newcommand{\n}{\lie{n}}
\newcommand{\m}{\lie{m}}

\newcommand{\AAA}{{\mathsf a}}
\newcommand{\Ul}{{\mathsf U}}
\newcommand{\ul}{{\mathsf u}}
\newcommand{\zl}{{\mathsf z }}

\newcommand{\fH}{{\mathfrak{H}}}
\newcommand{\fU}{{\mathfrak{U}}}

\newcommand{\A}{ \mathcal{A} }
\newcommand{\cT}{ \mathcal{T} }
\newcommand{\cE}{ \mathcal{E} }
\newcommand{\F}{ \mathcal{F} }
\newcommand{\G}{ \mathcal{G} }
\renewcommand{\O}{ \mathcal{O} }
\newcommand{\Ohat}{\widehat{ \mathcal{O} }}
\newcommand{\Nt}{{ \widetilde{\mathcal{N}} }}
\newcommand{\gt}{{ \widetilde{\mathfrak{g}} }}
\newcommand{\bc}{\mathsf{block}}
\newcommand{\iso}{{\;\stackrel{_\sim}{\longrightarrow}\;}}

\newcommand\conv{\star}
\newcommand\convf{\overset{!}{\star}}

\newcommand\al{\alpha}
\newcommand\la{\lambda}

\newcommand\bu{\bullet}

\newcommand{\Hom}{{\operatorname{Hom}}}
\newcommand{\Rep}{{\operatorname{Rep}}}
\newcommand{\uuHom}{{\underline{\operatorname{Hom}}}}
\newcommand{\Res}{{\operatorname{Res}}}
\newcommand{\Ind}{{\overline{\operatorname{Ind}}}}
\newcommand{\Mod}{{\operatorname{Mod}}}
\newcommand{\bimod}{{\operatorname{bimod}}}
\newcommand{\opmod}{{\operatorname{mod}}}

\newcommand{\oplusl}{\bigoplus\limits}

\newcommand{\imbed}{\hookrightarrow}

\newcommand{\Lotimes}{\overset{\operatorname{L}}{\otimes}}

\begin{document}

\title[]
{The small quantum group and the Springer resolution}

\author{R.Bezrukavnikov }

\address{MIT, Department of Mathematics, 77 Massachusetts ave., Cambridge MA 
02139, USA}

\email{bezrukav\@math.mit.edu}

\author{A.Lachowska}

\address{\'Ecole Polytechinique F\'ed\'erale de Lausanne, SB IMB CAG, 
B\^atiment MA, Station 8, Lausanne 1015, Switzerland }

\email{anna.lyakhovskaya\@epfl.ch}

\date{\today}

\begin{abstract}
In \cite{ABG} the derived category of the principal block in modules
over
the Lusztig quantum algebra at a root of unity is related to the
derived category of equivariant coherent sheaves on the Springer
resolution $\Nt$. In the present paper we deduce a similar relation
between the derived category of the principal block for the {\em small
(reduced) quantum algebra} $\ul$ and the derived category of
(non-equivariant) coherent sheaves on $\Nt$. As an application we get
 a geometric
description of Hochschild cohomology (in particular, the center) of
the regular block for $\ul$, and use it to give an explicit
description of a certain subalgebra in the center (such a  subalgebra 
was obtained
previously by another method and under more restrictive assumptions
in \cite{La}). We also briefly explain the relation of our result to
the geometric description \cite{BK} of the derived category of
modules over the De Concini -- Kac quantum algebra.
\end{abstract}


\maketitle

\centerline{\em To the memory of Iosef Donin.}

\vskip 0.3cm

\tableofcontents


\section{Quantum algebras}
\subsection{Basic notations and conventions}

Throughout the paper $\k$ is an algebraically closed field of
zero characteristic.

Let $R$ be a finite reduced root system in a 
$\k$-vector space $E$ and  fix
a basis of simple roots $S =\{\al_i, i \in I\}$. Let $\check\al$ denote
the coroot corresponding to the root $\alpha \in R$. The Cartan matrix
is given by $a_{ij} = \langle \al_i, \check\al_j \rangle$, where
$\langle \cdot , \cdot \rangle$ is the canonical pairing
$E^* \times E \to \k$. Let $W$ be the Weyl group of $R$.
 There exists a unique
$W$-invariant scalar product in $E$ such that $(\al, \al)=2$ for
any short root $\al \in R$. Set $d_i = \frac{1}{2} (\al_i, \al_i)
\in \{ 1,2,3 \}$ for each $i \in I$.
We denote by $\Y =\Z R$ the root lattice, and by
$\X = \{ \mu \in E : \langle \mu, \check\al \rangle 
\in \Z\;\; \forall \; \al \in R \}$
the weight lattice corresponding to $R$. The coweight lattice is
$\check\Y = \op{Hom}(\Y, \Z) \in E^*$. Let $R^+$ be the set of positive roots,
define the dominant weights by
$\X^+ = \{ \mu \in \X : \langle \mu, \check\al \rangle
\geq 0 \;\; \forall \; \al \in R^+ \}$
and set $\Y^+ = \Y \cap \X^+$.

Let $G$ be a connected semisimple group of adjoint type over $\k$
with the Lie algebra $\g$ corresponding to the root system $R$.
Let 
 $B$ be a Borel subgroup in $G$,
 and $N$ its unipotent radical.
Let  
 $\b$ and $\n$ be their respective Lie algebras.

\subsection{Quantum algebras at a root of unity}

Let $\k(q)$ denote the field of rational functions in the variable $q$.
We denote by $U_q(\g) = U_q $ the Drinfeld-Jimbo quantized enveloping algebra
of $\g$. It is generated over $\k(q)$ by $E_i, F_i, i \in I$ and
$K^{\pm 1}_\mu, \mu \in \check\Y$ subject to well-known relations,
see e.g. \cite{L}. We will write $K_i$ for $K_{d_i \check \al_i}$.
The algebra $U_q$ is a Hopf algebra over $\k(q)$.

Fix an odd positive integer $l$ which is greater than the Coxeter
number of the root system,
prime to the index of connection $|\X/\Y|$ and
prime to $3$ if $R$ has a component
of type $G_2$. Choose a primitive $l$-th root of unity $\xi \in \k$
and let $\A \subset \k(q)$ be the ring localized at $\xi$, and $\m$
the maximal ideal of $\A$. For any $n \in \N$ set
$[n]_d = \frac{q^{dn} - q^{-dn}}{q^d -q^{-d}}$ and
$[n]_d ! = \prod_{s=1}^n \frac{q^{ds} - q^{-ds}}{q^d -q^{-d}}$.

In $U_q$ consider the divided powers of the generators
$E_i^{(n)} = E_i^n /[n]_{d_i}!, F_i^{(n)} = F_i^n/[n]_{d_i}!,
i \in I, n \geq 1$,
and $\big[{ {K_\mu,m} \atop{n}} \big] $ as defined in \cite{L}. The
Lusztig's integral form $\Ul_\A$ is defined as an $\A$-subalgebra of
$U_q(\g)$ generated by these elements. $\Ul_\A$ is a Hopf subalgebra of
$U_q$. The Lusztig quantum algebra
at a root of unity $\Ul$ is defined by specialization of $\Ul_\A$ at $\xi$:
 $\Ul = \Ul_\A / \m \Ul_\A$. It has a Hopf algebra structure over $\k$.

Another version of the quantum algebra $\fU_\A$ was introduced in \cite{DK}.
This is the $\A$-subalgebra of $U_q$ generated by $E_i, F_i,
\frac{K_i -K_i}{q^{d_i} -q^{-d_i}}, \; i \in I$ and $K_\mu,
\mu \in \check\Y$. The elements $K_i^l$ are central in $\fU_\A$.
The De Concini-Kac quantum algebra is defined as
$\fU = \fU_\A/(\m \fU_\A + \sum_{i \in I}(K_i^l -1)\cdot \fU_\A)$.
It is a Hopf algebra over $\k$.

By definition $\fU_\A \subset \Ul_\A$. After the specialization at $\xi$, 
the imbedding 
of $\A$-forms induces a (not injective)
Hopf algebra homomorphism $\fU \map \Ul$. The image of this homomorphism
is the {\it small quantum group} $\ul$. Equivalently, $\ul$ is a subalgebra
in $\Ul$ generated by the elements  $E_i, F_i,
\frac{K_i -K_i}{q^{d_i} -q^{-d_i}}, \; i \in I$ and $K_\mu,
\mu \in \check\Y$. Since we have assumed $l$ to be odd, $\ul$
is a Hopf algebra over $\k$.

\subsection{Quantum Lusztig-Frobenius map}

Let $\{e_i, f_i, h_i \}_{i \in I}$ be the standard Chevalley generators
of the Lie algebra $\g$. Lusztig proved that the map
$E_i^{(l)} \map e_i$, $E_i \map 0$, $F_i^{(l)} \map f_i$, $F_i \map 0$
for all $i \in I$ can be extended to a well-defined surjective
algebra homomorphism
$\phi: \Ul \map \hat{U}(\g)$ which is called the {\it quantum Frobenius map}.
Here $\hat{U}(\g)$ stands for a certain completion of the universal enveloping
$U(\g)$ (see \cite{L}) such that the representation category of finite
dimensional $\hat{U}(\g)$-modules may be identified with that of the group
$G$.
The kernel of this map coincides with the two-sided ideal
in $\Ul$ generated by the augmentation ideal of the small quantum group,
$\ul_\epsilon$. One has an exact sequence of algebras

\be{Fro} 0 \map (\ul_\epsilon) \map \Ul \stackrel{\phi}\map \hat{U}(\g) .
\end{equation}

We let $\Rep(G)$ denote the category of finite dimensional algebraic
$G$-modules, and $\Rep(\Ul)$ be the category of finite dimensional
$\Y$-graded
$\Ul$-modules.

The pull-back of the Frobenius homomorphism gives rise to the
functor between the tensor categories:

\be{Ffro} \phi^* : \op{Rep}(G) \map \op{Rep}(\Ul), \;\;\; V \map V^\phi.
\end{equation}

\section{Functor to a derived category of  $\ul$-modules }

In the next three subsections we recall the construction of the functor
introduced
in \cite{ABG} from a derived  category of coherent $G$-equivariant sheaves
on the Springer resolution corresponding to $G$ to a certain derived
category of representations of $\Ul$.

\subsection{The principal block} \label{blockU}

The category $\op{Rep}(\Ul)$  is an abelian
artinian category, and therefore is a direct sum of its indecomposable
abelian subcategories, or blocks.

 We write $L(\nu)$ for a simple finite dimensional $\Ul$-module of
highest weight $\nu \in \Y^+$. For any $\la \in \Y$ let
$L_\la $ be the finite dimensional simple $\Ul$-module with
highest weight $w(l\la +\rho) -\rho$, where $\rho = \frac{1}{2}
\sum_{\al \in R^+} \al$ and $w \in W$ is the unique element such
that $w(l\la +\rho) -\rho$ is a dominant weight.

We write $\bc (\Ul)$ for the block of $\op{Rep}(\Ul)$ which contains
the trivial representation (the principal block). Equivalently, $\bc(\Ul)$
is the full subcategory of the abelian category of left $\Ul$-modules
formed by the modules $M$ such that all simple subquotients of $M$
are of the form $L_\la$, $\la \in \Y$.

We let $D^b \bc(\Ul)$ denote the corresponding bounded derived category.

\subsection{Springer resolution and $\Ul$-modules}
Using the adjoint action of $B$ on the nilradical $\n$ of $\b
=\op{Lie}(B)$, one defines the {\it Springer resolution} $\Nt = G
\times_{B} \n$ as the quotient of $G \times \n$ by the action $h
\cdot(g,x) = (g h^{-1}, Ad_h(x))$. Thus $\Nt$ is an
algebraic variety equipped with an algebraic action of $G$.
 The multiplicative group $\C^*$  acts on
$\Nt$ by $t:(g,x)\mapsto (g,t^2x)$ along the fibers.

 Let $Coh^{G \times \C^*}(\Nt)$ (resp.
$Coh^{\C^*}(\Nt)$) denote the abelian category of $G \times
\C^*$-equivariant (respectively $\C^*$-equivariant) coherent sheaves on
$\Nt$. For an abelian category $C$ we write $D^b C$ for its bounded
derived category.

The following result was obtained in \cite{ABG}, Corollary 1.4.4.

\begin{theorem} \label{F}
There exists a triangulated functor
$$
F: D^bCoh^{G \times \C^*} (\Nt) \longrightarrow D^b\bc(\Ul),
$$
such that:
\begin{enumerate}
\item{$ F(z^i \otimes \mathscr{F}) = F(\mathscr{F})[i]$ for any
${\mathscr F} \in D^b Coh^{G \times \C^*}(\Nt)$ and any $i \in \Z$.}
\item{The functor $F$
induces, for any ${\mathscr{F}},{\mathscr{F}}'
\in D^b Coh^{G\times\C^*}(\Nt),$ canonical isomorphisms
$$\bigoplus_{i\in\Z}\,\op{Hom}^{\bullet}_{D^b Coh^{G \times \C^*}(\Nt)}
({\mathscr{F}},z^i\otimes{\mathscr{F}}') \iso
\op{Hom}^{\bullet}_{D^b \bc(\Ul)}(F({\mathscr{F}}),
F({\mathscr{F}}')).$$}
\item{The image of $F$ generates the target category as a
triangulated category.}
\item The functor intertwines the natural action of the tensor
category $\op{Rep}(G)$ on $D^bCoh^{G\times \C^*}(\Nt)$ with the
action on $D^b\bc(\Ul)$ coming from the Lusztig-Frobenius
homomorphism $\phi$; i.e., for every $V\in \Rep(G)$ we have an
isomorphism $F(V\otimes \F)\cong \phi^*(V)\otimes F(\F)$, satisfying
the natural compatibilities.
\end{enumerate}
\end{theorem}

Here $z^i \otimes \mathscr{F}$ denotes the $\C^*$-equivariant sheaf
$\mathscr F$ with the $\C^*$-equivariant structure twisted by the
character $z  \to z^i$. The notation $\mathscr{ F}[k]$ stands for
the homological shift of ${\mathscr F}$ by $k$ in the derived category.

\subsection{Springer resolution and $\ul$-modules}
We would like to obtain a similar functor into a derived category of
$\ul$-modules.

We let $\Rep(\ul)$ be the category of finite dimensional
$\ul$-modules. The category $\op{Rep}(\ul)$ is an abelian artinian
category, and therefore is a direct sum of blocks. Let $\bc(\ul)$ 
denote the block containing the trivial representation.

 \begin{theorem} \label{F_u}
There exists a triangulated functor
$$ F_u : D^bCoh^{\C^*} (\Nt) \longrightarrow D^b\bc(\ul) $$
fitting into the commutative diagram
$$ \begin{array}{ccc}
 D^b Coh^{G\times\C^*}(\Nt)& \stackrel{F}{\longrightarrow} &
D^b\bc(\Ul)
\\
                           &                               &
  \\
 \op{Forget} \downarrow   &                  & \downarrow \op{Res} \\
                           &                         &
  \\
D^b Coh^{\C^*}(\Nt)&\stackrel{F_u} {\longrightarrow} & D^b\bc(\ul)  \\
       \end{array} $$
The functor $F_u$ satisfies properties (1--3) stated in Theorem
\ref{F}.
\end{theorem}


\proof We need to recall a relation between the categories $\bc(\Ul)$ and
$\bc(\ul)$ established in  \cite{AG}.\footnote{\cite{AG} concentrate
 on the case
 of an {\em even} (and divisible by 3 if $G$ has a component of type $G_2$)
root of
unity. However, the simpler case of an odd root of unity
considered in the present paper is also covered by the general
Theorem 2.8 of {\em loc. cit.} (it is easy to show that the
assumptions of \cite{AG}, Theorem 2.8 are satisfied in our case).}

Recall that we
have the restriction functor $\Res:\Rep(\Ul)\to \Rep(\ul)$. The right
adjoint functor to $\Res$ is defined if we pass to the categories of
ind-objects. The category on ind-objects in $\Rep(\Ul)$ can be identified
with the category of locally finite $\Y$-graded
modules over $\Ul$; we denote it by $\Rep^{lf}(\Ul)$.
 Thus we have the ``locally finite induction'' functor $\Ind:\Rep(\ul)
\to \Rep^{lf}(\Ul)$.
An object in $\Rep^{lf}(\Ul)$
 $\Ind(M)$, $M\in \bc(\ul)$ carries an
additional structure; namely, the object $\Ind(M)$ of the tensor
category $\Rep^{lf}(\Ul)$ is naturally a module over
 the algebra $\O(G)$.
 Here $\O(G)$ is the
algebra of regular functions on the algebraic group $G$ viewed as a
algebra in the tensor category of $\g$-modules. To make sense of
the notion of an $\O(G)$-module in $\Rep^{lf}(\Ul)$ we need to fix
an action of the tensor category $\Rep(G)$  on  $\Rep^{lf}(\Ul)$.
Such an action is given by $V:M\mapsto V^\phi \otimes M$.

A module structure on a locally finite module $M$
amounts to a collection $h_V$  of isomorphisms
 $V^\phi \otimes M \to \underline{V} \otimes M$ fixed for every
algebraic $G$-module $V$; here $\underline{V}$ denotes the vector
space underlying the representation $V$. The collection of
isomorphisms $h_V$ has to satisfy a certain compatibility condition
spelled out, e.g., in \cite{AG}.

For a finite dimensional module $M\in \Rep(\ul)$ the
 $\O(G)$ module $\Ind(M)$ is finitely generated.
 We let $\Mod_{\Rep^{lf}(\Ul)}(\O(G))$ denote
the category of finitely generated $\O(G)$ modules in
$\Rep^{lf}(\Ul)$.

It follows from the results of \cite{AG}
(or from the much more
general Barr-Beck theorem, see e.g. \cite{Ml})
 that the functor from $\Rep(\ul)$
to  $\Mod_{\Rep^{lf}(\Ul)}(\O(G))$ is an
equivalence.

It is not hard to show that the functor $\Res$ sends 
$\bc(\Ul)$ to $\bc(\ul)$, while $\Ind$ sends
 of $\bc(\ul)$ to $\bc^{lf}(\Ul)$, where $\bc^{lf}(\Ul)$ is the category
of ind-objects in $\bc(\Ul)$ (identified with the category of locally finite
graded $\Ul$-modules, which are unions of modules in $\bc(\Ul)$).
It follows that we have an equivalence
$$
\bc(\ul)\cong \Mod_{\bc^{lf}(\Ul)}(\O(G)),$$
where $\Mod_{\bc^{lf}(\Ul)}(\O(G))$ is the category of
finitely generated $\O(G)$-modules in the category $\bc^{lf}(\Ul)$.
(Notice that the action
$\Rep(G)$ on $\Rep^{lf}(\Ul)$ preserves
 $\bc^{lf}(\Ul)$, thus  the notion of an $\O(G)$-module in this category
is well-defined).

Furthermore,  $\Res$, $\Ind$ are exact, hence $\Ind$ sends injective
objects to injective ones. It is not hard to deduce that the above
equivalences are  inherited by the derived categories, i.e. we have
$$D^b\Rep(\ul)\cong
\Mod_{D^b\Rep^{lf}(\Ul)}(\O(G));$$
$$D^b\bc(\ul)\cong
\Mod_{D^b\bc^{lf}(\Ul)}(\O(G)).$$
Here $\Mod_{D^b\Rep^{lf}(\Ul)}(\O(G))$,
$\Mod_{D^b\bc^{lf}(\Ul)}(\O(G))$ denote the categories of 
finitely generated (equivalently, compact) $\O(G)$-modules
in the corresponding bounded derived categories.

A similar relation exists between the (derived) categories of
equivariant and non-equivariant (or equivariant under a smaller
group) coherent sheaves. More precisely, we have a pair of adjoint
functors
 $R:D^bCoh^{G\times \C^*}(\Nt)\to D^bCoh^{\C^*}(\Nt)$ and $Av:
 D^bCoh^{\C^*}(\Nt)\to D^bQCoh^{G\times \C^*}(\Nt)$, where $R$ stands
for the restriction of equivariance functor, and the right adjoint
$Av$ is the ``averaging'' functor $a_*pr^*$,
 where $pr,a:G\times \Nt\to \Nt$ are, respectively, the projection
 and the action map.
An object $Av(\F)$, $\F\in Coh^{\C^*}(\Nt)$ carries an additional
structure of a finitely generated $\O(G)$-module. As above, the
structure of an $\O(G)$-module on a quasi-coherent equivariant
sheaf amounts to the data of an isomorphism
 $h_V:V \otimes \F \iso \underline{V} \otimes \F$ fixed for every
algebraic $G$-module $V$ and satisfying the compatibilities of
\cite{AG}. Here the action of $\Rep(G)$ on $QCoh^{G\times
\C^*}(\Nt)$ is given by $V:\F \mapsto V\otimes \F$, where  the
$G$-equivariant structure on the sheaf $V\otimes \F$ is the tensor
product of the equivariant structure on $\F$ and the action of
$G$ on $V$.

It is elementary to check that this way we obtain an equivalence
$$Coh^{\C^*}(\Nt)\cong \Mod_{QCoh^{G\times \C^*}(\Nt)}(\O(G)),$$
where $\Mod_{QCoh^{G\times \C^*}(\Nt)}(\O(G))$ is the
category of finitely generated (compact) $\O(G)$-modules in $QCoh^{G\times
\C^*}(\Nt)$.

Moreover, the functors $\Res$, $\Ind$ are exact, and the functor
$\Ind$ sends injective ind-objects to injective ind-objects. It
follows that
$$D^bCoh^{\C^*}(\Nt)\cong \Mod_{D^bQCoh^{G\times \C^*}(\Nt)}
(\O(G)),$$ where the category in the right hand side is the category
of finitely generated (compact) modules for $\O(G)$ in $D^bQCoh^{G\times
\C^*}(\Nt)$.

It remains to notice that in view of property (4) of Theorem \ref{F}
the functor $F$ intertwines the actions of the tensor category
$\Rep(G)$, thus it induces a functor between the categories of
finitely generated modules. It is immediate to see that properties
(1--3) of the functor $F$ from Theorem \ref{F} yield similar
properties of the induced functor between the module categories.
$\square$

\subsection{Connection to the Kac-De Concini algebra
and a result of \cite{BK}} In this subsection we sketch an
alternative way to prove and somewhat strengthen Theorem \ref{F_u},
which relies on the result of \cite{BK} and elementary theory of
differential-graded schemes. The material of this subsection is not
used elsewhere in the paper; the details are omitted.

Recall that $\fU$ denotes the {\em  De Concini -- Kac} algebra.
The center $Z(\fU)$ contains two subalgebras, the {\em
Harish-Chandra center} $Z_{HC}$, and the {\em $l$-center} $Z_l$.
The Harish-Chandra center $Z_{HC}$ is obtained from the center
of the quantized enveloping algebra $U_q$ by specialization.
The center $Z(U_q)$ is isomorphic by the quantum
Harish-Chandra map to $((U_q)_0)^{W\ltimes \Gamma}$, where
$(U_q)_0$ is the subalgebra of $U_q$ generated by
$\{ K_\mu \}_{\mu \in \check\Y}$, and
$\Gamma$ is the group of homomorphisms $ \check\Y \to \{\pm 1 \}$.
The $l$-center $Z_l$ is the central subalgebra of $\fU_l$ generated
by the $l$-th powers of the generators $\{E_i, F_i\}_{i \in I}$,
$\{K_\mu\}_{\mu \in \check\Y}$. Then we have
$Z(\fU) \simeq Z_l \otimes_{Z_l \cap Z_{HC}} Z_{HC}$ (see e.g. \cite{DKP}).

Let
$\op{Rep}(\fU)_{0}$, respectively, $\op{Rep}(\fU)_{\widehat 0}$
denote the full subcategory in $\op{Rep}(\fU)$ consisting of
modules killed by the augmentation ideal in $Z_l$ (respectively, by
some power of this ideal). Similarly, let $\op{Rep}(\fU)^{0}$,
$\op{Rep}(\fU)^{\widehat 0}$ be the full subcategories in
$\op{Rep}(\fU)$ consisting of modules killed by the augmentation
ideal in $Z_{HC}$ (respectively, by some power of this ideal). We
also set $\op{Rep}(\fU)_{0}^{0}=\op{Rep}(\fU)_{0}\cap
\op{Rep}(\fU)^{0}$ etc.

 We have $\ul = \fU\otimes _{Z_l} \k$, thus
$\op{Rep}(\fU)_{0}=\op{Rep}(\ul)$.
  It is not hard to show that
  $\op{Rep}(\fU)_{0}^{\widehat 0}\cong \bc(\ul)$.

\medskip

 The main
theorem of \cite{BK} yields an equivalence \begin{equation}
\label{BKeq} D^b \op{Rep}(\fU)_{\widehat 0}^{\widehat 0}\cong D^b
Coh_{G/B}(\gt).\end{equation}
 Here $\gt =G\times ^B\b$, and $Coh_{G/B}(\gt)$ is the full
 subcategory in $Coh(\gt)$ consisting of sheaves set-theoretically
 supported on the zero section.
  The completion $\widehat{Z_l}$ of $Z_l$ at the
 augmentation ideal is naturally identified with the completion $\Ohat(\g)$
 of
 the polynomial algebra $\op{Sym}(\g)$ at the augmentation ideal
 (here we use an identification $\g\cong \g^*$ provided by an invariant
 quadratic form).
 The equivalence \eqref{BKeq} intertwines the action of
 $\widehat{Z_l}$ on the derived $\fU$-module category with the
 action of $\Ohat(\g)$ on the derived category of coherent sheaves
 coming from the Grothendieck-Springer map $\gt\to \g$.

One can show that the completion of $\fU$ at the augmentation
ideal of $Z_l\cdot Z_{HC}$ is flat over $\widehat{Z_l}$.

Furthermore, one can deduce by 
a base change argument an equivalence
\begin{equation}
\label{derpr} D^b \op{Rep}(\fU)_{0}^{\widehat 0}\cong DGCoh
(\gt\overset{\op{L}}{\times}_\g \{0\}).\end{equation}

Here $\gt\overset{\op{L}}{\times}_\g \{0\}$ is the {\em differential
graded (DG) scheme}, which is the {\em derived fiber product} of the
schemes $\gt$ and $\{0\}=Spec(\k)$ over $\g$, while $DGCoh$ denotes
the derived category of sheaves of DG  $\O$-modules over the DG
scheme (see e.g. \cite{Ka} for the definitions).

By the definition, the structure sheaf of the DG scheme
$\gt\overset{\op{L}}{\times}_\g \{0\}$ is a sheaf of DG-algebras on
$G/B$, which is well-defined up to a quasiisomorphism. A possible
construction of a representative of the quasi-isomorphism class is
as follows: $\O_{\gt\overset{\op{L}}{\times}_\g \{0\}}=\pi^*{\mathbb
K}_\g$, where $\pi:\gt\to \g$ is the Grothendieck-Springer map, and
${\mathbb K}_\g=\Lambda^\bu(\g)\otimes \op{Sym}^\bu(\g)$ is the
Koszul complex of $\g$.

It is not hard to show that for a vector bundle (locally free sheaf)
$\cE$ on an algebraic variety together with an embedding of vector
bundles $\cE\subset V\otimes \O$ of $\cE$ into the trivial vector
bundle, the sheaf of DG-algebras
$\O_{E\overset{\op{L}}{\times}_V\{0\}}$ is canonically
quasiisomorphic to the sheaf of DG-algebras with zero differential
$$\Lambda ((V\otimes
\O/\cE)^*[1])=\underline{\op{Tor}}_{\bu}^{\O(V)}(\O_E,\k).$$ Here
$E$ denotes the total space of $\cE$.

In particular we see that $\O_{\gt\overset{\op{L}}{\times}_\g
\{0\}}$ can be represented by the DG algebra with zero differential
$
\Lambda^\bu(\Omega^1_{G/B}[1])$, where $\Omega^1_{G/B}$ is the
locally free sheaf of 1-forms on $G/B$. 
Thus we get an
equivalence $$D^b\bc(\ul)\cong
DGCoh(\Lambda^\bu(\Omega^1_{G/B}[1])).$$

Finally, the standard Koszul (or $S-\Lambda$) duality, see e.g.
\cite{BGS} (cf. also \cite{ABG}, \S 3.3), gives a canonical equivalence
$$DGCoh(\Lambda^\bu(\Omega^1_{G/B}[1]))\cong
DGCoh(\op{Sym}^\bu(T_{G/B}[-2])),$$ where $T_{G/B}$ is the
tangent sheaf of $G/B$. Thus we get
$$D^b\bc(\ul)\cong
DGCoh(\op{Sym}^\bu(T_{G/B}[-2])).$$ Notice that the relative
spectrum of the sheaf of commutative rings $\op{Sym}^\bu(T_{G/B})$
on $G/B$ is nothing but $\Nt$. Thus the last equivalence implies
Theorem \ref{F_u}.

This method can also be used to provide a similar description for
the derived category of the regular block for the restricted
enveloping algebra of a semi-simple Lie algebra over a field of
positive characteristic; the reference to \cite{BK} should then be
replaced by a reference to \cite{BMR}.

\section{Hochschild cohomology of the principal block of $\ul$}


The finite dimensional Hopf algebra $\ul$ decomposes as a left $\ul$-module
into a finite direct sum of finite dimensional submodules.
Denote by $\ul_0$ the largest direct summand for which all its
simple subquotients belong to the principal block of the category
$\op{Rep}(\ul)$. Then $\ul_0$ is a two-sided ideal in $\ul$,
which will be called the principal block of $\ul$.

Let $\zl$ denote the center of $\ul$. It decomposes into a direct
sum of ideals according to the block decomposition of $\ul$.
Set $\zl_0 = \zl \cap \ul_0$.

The rest of this section contains 
a computation of the Hochschild cohomology of the principal block of
$\ul$ and a description of the center $\zl_0$.

\subsection{The result}

Recall that we have a $G$-equivariant isomorphism of vector bundles
$G \otimes_{B} \n \cong T^*(G/B) = \Nt$ and that the
multiplicative group acts on $\Nt$ by dilations on the fibers: an
element $t \in \C^*$ acts on $n$ by multiplication by $t^2$.
Consider the  coherent sheaf of poly-vector fields $\Lambda^\bullet
T(\Nt)$ on $\Nt$. The direct image of this sheaf to $G/B$ is in fact
bi-graded. Here the first grading is the natural grading
$\Lambda^\bullet T(\Nt) =\oplus_{j = 0}^{\op{dim}(\Nt)} \Lambda^j
T(\Nt)$. The second grading comes from the induced action of $\C^*$
on $\Nt$. We will write $\Lambda^j T(\Nt)^k$ for the $(j,k)$-th
component with respect to this bi-grading; this is a locally free
$G$-equivariant coherent sheaf on $G/B$. Notice that
$\Lambda^j T(\Nt)^k=0$ for odd $k$.

\begin{theorem} \label{Hoh}
There exists an isomorphism of algebras between
the total Hochschild cohomology of the principal block $\ul_0$ and
the total cohomology of $\Nt$ with coefficients
in $\Lambda^\bullet T(\Nt)$; here the algebra structure on the second space
comes from multiplication in the exterior algebra $\Lambda^\bullet T(\Nt)$.
The isomorphism is compatible with the grading as follows:
$$  \op{HH}^s (\ul_0) \cong \oplus_{i+j+
k=s} \op{H}^i(\Lambda^j T(\Nt))^k. $$
\end{theorem}

\begin{remark}
``Morally'' the Theorem is a direct consequence of Theorem \ref{F_u}.
More precisely, suppose it were possible to define for a triangulated category
$\cT$ a triangulated category of endo-functors $End(\cT)$ with good properties.
Some of the expected properties are as follows: if $\cT$ is the derived
category of modules over an algebra, then $End(\cT)$ is the derived category
of bi-modules over the same algebra; while if $\cT$ is the  derived category
of (equivariant) coherent sheaves on an algebraic variety $X$, then
$End(\cT)$ is the derived category of (equivariant) coherent sheaves on $X^2$.
The relation between the categories $D^b(Coh^{\C^*}(\Nt))$ and $D^b(\bc(\ul))$
explained in Theorem \ref{F_u} would then imply a similar relation
between the endomorphism categories. Expanding  property (2) of Theorem
\ref{F} for $\F$ and $\F'$ being the identity functor we would get Theorem
\ref{Hoh}.

It is well known that the naive category of endo-functors of a triangulated
category does not, in fact, carry a natural triangulated structure and does
not satisfy the properties indicated above. One can probably define an
appropriate category of endomorphisms by working with differential graded
categories (or in another rigid setting, such as that of
$A_\infty$ categories). We found it more effective to derive Theorem
\ref{Hoh} by a more elementary ad hoc argument.
\end{remark}

\begin{corollary} \label{zl0} \label{center}
The principal block of the center of $\ul$ is isomorphic as an algebra to
$$ \zl_0 \cong \oplus_{i+j+
k=0} \op{H}^i(\Lambda^j T(\Nt))^k. $$
\end{corollary}

The Corollary is immediate by setting $s=0$ in Theorem \ref{Hoh}.

The proof of  Theorem \ref{Hoh} is based on the following standard
statement, which is 
an algebraic version of the Hochschild-Kostant-Rosenberg Theorem, see,
e.g. \cite{Swan}.
\begin{lemma} 
 Let $X$ be a smooth variety and $\delta: X \to X \times X$ the diagonal 
imbedding. Then there is an algebra isomorphism 
$$ \op{Ext}_{Coh(X \times X)}^q (\delta_* \mathcal{O}_X,
\delta_* \mathcal{O}_X) \cong \bigoplus_{i+j=q} H^i(\Lambda^j T_X). $$
\end{lemma}

In what follows we identify $\ul\cong \ul^{op}$ by means of the
antipode of the Hopf algebra $\ul$, and we identify the category of
$\ul$-bimodules with that of $\ul\otimes \ul$-modules. We let
$\bc(\ul^2)$ denote the block of the trivial representation in the
category of $\ul$-bimodules.

 In view of the Lemma, Theorem \ref{Hoh} follows immediately
from the following

\begin{proposition}\label{Phi}
There exists a functor $\Phi: D^b Coh^{\C^*}(\Nt \times \Nt)
\rightarrow D^b\bc(\ul^{\otimes 2})$, which satisfies the properties
(1--3) of Theorem \ref{F}, and such that
$\Phi(\delta_*(\O_\Nt))\cong R$; here $R$ is the maximal summand in
the regular bimodule for $\ul$ belonging to the principal block, and
$\delta:\Nt\to \Nt\times \Nt$ is the diagonal embedding.
\end{proposition}

\subsection{The proof}
The rest of the section is devoted to the proof of Proposition
\ref{Phi}. We start with some auxiliary statements.

Recall a   monoidal structure on the derived category of bimodules
and an action of this monoidal category on the derived category of
modules. More precisely, let $\AAA$ be an associative ring, and
$\AAA-\op{mod}$, $\AAA-\op{bimod}$ be the categories of
$\AAA$-modules and of $\AAA$-bimodules respectively. We set $B\conv
M=B\Lotimes_\AAA M$; here $B\in D^-(\AAA-\bimod)$,  and $M$  is
either an object of the same category, or $M\in D^-(\AAA-\opmod)$.

In the first case we get a monoidal structure on $D^-(\AAA-\bimod)$,
while the second one gives an action of this monoidal category on
$D^-(\AAA-\opmod)$.

We also have the dual operation $B\convf M=R\Hom_\AAA(B,M)$. This
formula defines functors $D^+(\AAA-\bimod)\times D^+(\AAA-\bimod)\to
D^+(\AAA-\bimod)$, $D^+(\AAA-\bimod)\times D^+(\AAA-\opmod)\to
D^+(\AAA-\opmod)$.

\begin{lemma}\label{AAA} Let $\AAA$ be a (left and right) Noetherian
 associative ring, and $B\in D^b(\AAA-\bimod)$ be such that the image of
 $B$ in $D^b(\AAA-\opmod)$, $D^b(\AAA^{op}-\opmod)$ under the
 functors of forgetting the right (respectively, left) action is a
 perfect complex (i.e. can be represented by a finite complex of
 finitely generated projective modules).

a) The functor $M\mapsto B\conv M$ from $D^b(\AAA-\opmod)$ to
$D^b(\AAA-\opmod)$ has a right adjoint given by $M\mapsto B\convf M$.

b) We have a canonical isomorphism $$(B\convf C)\conv D\cong B\convf
(C\conv D)$$
Here  $C\in D^b(\AAA-\bimod)$,
and $D$ either lies in $D^-(\AAA-\bimod)$, or in $D^-(\AAA-\opmod)$.

c)  
 Assume moreover that the functor
$D^b(\AAA-\opmod)\to D^b(\AAA-\opmod)$,
 $M\mapsto B\conv M$ is an equivalence.
Then the functor $D^b(\AAA-\bimod)\to D^b(\AAA-\bimod)$, $M\mapsto
B\convf M$, is an equivalence sending $B$ to the regular bimodule.
\end{lemma}

\proof a) and b) are standard. To check (c) observe that right
adjoint to an equivalence is the inverse equivalence. Thus the
composition of endo-functors $M\mapsto B\conv M$ and $M\mapsto
B\convf M$ of $D^b(\AAA-\opmod)$
is isomorphic to identity. This composition is given by
$$M\mapsto B\convf (B\conv M)\cong (B\convf B)\conv M,$$
where the isomorphism is provided by part (b).

Thus setting $C=(B\convf B)$ we see that the endo-functor of
$D^b(\AAA-\opmod)$, $M\mapsto C\conv M$ is isomorphic to identity.
This is easily seen to imply that $C$ is isomorphic to the regular
bimodule. It remains to show that the endofunctor of
$D^b(\AAA-\bimod)$, $M\mapsto B\convf M$ is an equivalence.

Its left adjoint functor $M\mapsto B\conv M$ (where the adjunction is provided
by part (a) of the Lemma) is also its right inverse: this is clear from (b)
and the established isomorphism between $B\convf B$ and the regular bimodule.
This implies that the endo-functor of $D^b(\AAA-\bimod)$, $M\mapsto B\conv M$
is a full embedding, and the category  $D^b(\AAA-\bimod)$ admits
a semi-orthogonal decomposition $$D^b(\AAA-\bimod)=\op{Im}(M\mapsto B\conv M)
*\op{Im}(M\mapsto B\conv M)^\perp.$$
Here $\op{Im}(M\mapsto B\conv M)^\perp\subset D^b(\AAA-\bimod)$
 is the full subcategory
of objects $N\in D^b(\AAA-\bimod)$ satisfying $\Hom(B\conv M,N)=0$
for all $M$, and the displayed formula means that
 for every object $M\in D^b(\AAA-\bimod)$
there exists a (canonical) distinguished triangle
$$M_1\to M\to M_2, \ \ \ \ \
M_1\in \op{Im}(M\mapsto B\conv M),\ M_2\in \op{Im}(M\mapsto B\conv
M)^\perp.$$ Thus to show that $M\mapsto B\conv M$, $M\mapsto B\convf
M$ are auto-equivalences of $D^b(\AAA-\bimod)$ it is enough to show
$\op{Im}(M\mapsto B\conv M)^\perp=\{0\}.$ We have
 $\op{Im}(M\mapsto B\conv M)^\perp=\op{Ker}(M\mapsto B\convf M).$
Thus it suffices to show that $B\convf M\ne 0$ provided $M\ne 0$.
The endo-functors of $D^b(\AAA-\bimod)$, $D^b(\AAA-\opmod)$ defined
by $M\mapsto B\convf M$ are obviously compatible via the functor of
forgetting the right action. Since the endo-functor of
$D^b(\AAA-\opmod)$, $M\mapsto B\convf M$ is an auto-equivalence, it
kills no non-zero objects; hence the same is true for the
endo-functor of $D^b(\AAA-\bimod)$. $\square$

\medskip

Replacing $\g$ by $\g \oplus \g$ in Theorem \ref{F_u}, we obtain the
functor
$$ F_{u^2} : D^b Coh^{\C^*}(\Nt \times \Nt)
\rightarrow D^b\bc(\ul)^{\otimes 2}$$ satisfying properties (1--3)
of Theorem \ref{F}.

\begin{lemma}\label{Bgut}
Set $B=F_{u^2}(\delta_*(\O_\Nt))$. Then we have a functorial
isomorphism $$F_u(\F)\cong B\conv  F_u(\F\check{\ })^*$$ for $\F\in
D^bCoh^{\C^*}(\Nt)$; here $\F\check{\ }=R\uuHom (\F,\O_{\Nt})$ is
the dual sheaf, and $*$ denotes the usual duality for (complexes of)
modules.
\end{lemma}

\proof It will be convenient to use the abbreviation $\Hom^\bu(M,N)=
\oplusl_i \Hom(M,N[i])$ for $M,N\in D^b\bc(\ul)$ and
$\Hom^\bu(\F,\G)= \oplusl_{i,j} \Hom(\F,z^j\otimes \G[i])$ for
$\F,\G\in D^bCoh^{\C^*}(\Nt)$. Thus property (2) of Theorem \ref{F}
asserts that $$\Hom^\bu(\F,\G)\iso \Hom^\bu(F_u(\F),F_u(\G))$$ for
$\F,\G\in D^bCoh^{\C^*}(\Nt)$. Notice that this isomorphism
preserves the grading, where elements of $\Hom(M,N[i])$ are assigned
degree $i$, while elements of $\Hom(\F, z^j\otimes \G[i])$ are assigned
degree $i+j$.

It suffices to construct a functorial isomorphism
$$ \Hom^\bu(\F,\G)\cong \Hom^\bu(B\conv
 F_u(\F\check{\ })^*, F_u(\G)),
 $$
 preserving the grading.
 Indeed, then plugging in $\F=\G$ we get a morphism
 $F_u(\F\check{\ })^*\to F_u(\F)$, which induces an isomorphism
 between the spaces of homomorphisms to $F_u(\G)$ for any $\G$.
Since the image of $F_u$ generates $D^b\bc(\ul)$, we see that this
morphism is an isomorphism.

We have:
\begin{multline*}
\Hom^\bu((B\conv
 F_u(\F\check{\ })^*, F_u(\G))\cong \Hom^\bu_{D^b\bc(\ul^{2})}(B,
 F_u(\F\check{\ })\boxtimes F_u(\G))\cong \\
 \Hom^\bu(\delta_*(\O_\Nt), \F\check{\ }\boxtimes \G);
 \end{multline*}
here we used the obvious compatibility of $F_u$ with external
products, and also the adjunction $\Hom^\bu(B\conv M, N)\cong
\Hom^\bu(B, M^*\boxtimes N)$ valid for any bimodule $B$ and modules
$M,\, N$ (or complexes of such). Applying duality $\F\mapsto R\uuHom(\F,
\O_{\Nt^2})$ we get  $$\Hom^\bu(\delta_*(\O_\Nt), \F\check{\
}\boxtimes \G)\cong \Hom^\bu(\F\boxtimes \G\check{\ },
\delta_*(\O_\Nt)\check{\ }).$$ An easy calculation shows that $
\delta_*(\O_\Nt)\check{\ }\cong z^{2d}\otimes
\delta_*(\O_\Nt)[-2d]$, where $d=\dim G/B$. Thus the latter Hom
space can be rewritten as
$$\Hom^\bu(\F\boxtimes \G\check{\ },
z^{2d}\otimes \delta_*(\O_\Nt)[-2d])\cong \Hom^\bu(\F, z^{2d}\otimes
\G[-2d])\cong \Hom^\bu(\F,\G).$$ It is clear that the last
isomorphism preserves the grading (the twists of the $\C^*$ action
and the homological shifts cancel). $\square$

{\it Proof of Proposition \ref{Phi}}. We start with the functor
$F_{u^2}$ (introduced before Lemma \ref{Bgut}), and set
$B=F_{u^2}(\delta_*(\O_\Nt))$ as in Lemma \ref{Bgut}.

 Lemma
\ref{Bgut} implies that $B$ is of finite projective dimension as
left and right $\ul$-module, and
 $M\mapsto B\conv M$ is an auto-equivalence of $D^b\bc(\ul)$.
 Indeed, objects of the form $F_u(\F\check{\ })^*$ generate
 $D^b\bc(\ul)$ as a triangulated category, since the image of
 $F_u$ does. Since the functor $M\mapsto B\conv M$ sends such
 objects into the bounded derived category (rather than into
 $D^-\bc(\ul)$), it also sends the whole of $D^b\bc(\ul)$ into
 itself. It follows that $B$ has a finite projective dimension
 as a right $\ul$-module. The involution of switching the two
 factors in $\ul\otimes \ul$ obviously sends $B$ into itself, thus
 $B$ also has a finite projective dimension as a left $\ul$-module.
Finally, the functor $M\mapsto B\conv M$ sends the set
$\{F_u(\F\check{\ })^*\ |\ \F\in D^b(Coh^{\C^*}(\Nt))\}$
  generating
 $D^b\bc(\ul)$ as a triangulated category into another generating
 set, and it induces an isomorphism on $\Hom$'s between objects in
 the generating set. Hence it is an equivalence.

Thus by Lemma \ref{AAA} there exists an autoequivalence $A$ of
$D^b\bc(\ul^2)$, sending $B$ to the regular bimodule. Then
$\Phi=A\circ F_{u^2}$ is readily seen to satisfy the required
properties. $\square$

\begin{remark}
One can ask for a more explicit description of the object
$B=F_{u^2}(\delta_*(\O_\Nt))$. This question is similar to the
question of describing the bimodule over the classical enveloping
algebra $U(\g)$ obtained as the global sections of the $D$-module
$\Delta_*(\O_{G/B})$; here $\Delta:G/B\to (G/B)^2$ is the diagonal
embedding, and $\Delta_*$ denotes the direct image in the category
of $D$-modules. In both cases one can show that the endo-functor of
the derived category of modules coming from this bimodule can be
described in terms of the action of the braid group on the derived
category of modules by {\em intertwining functors} \cite{BBcas}
(cf. also \cite{ABG}, \S 4.1). More precisely, it coincides with the action
of the canonical lifting to the braid group of the longest element in
the Weyl group. We neither prove nor use this fact in the present
paper.

\end{remark}
\subsection{An explicit subalgebra in the principal block of the center of the
small quantum group.}

We will describe an easily computable subalgebra in $\zl_0$. To
state the answer, define a commutative algebra $\fH$ of dimension
$2|W|-1$ as follows. Endow the space $H^\bu(G/B,\k)\oplus H^\bu(G/B,\k)$
with a commutative algebra structure given by $(h_1,h_2)\cdot
(h_1',h_2')=(h_1h_1', \epsilon (h_1)h_2'+\epsilon (h_1')h_2)$, where
$\epsilon:H^\bu(G/B)\to \k=H^0(G/B)$ is the augmentation. We let
$\fH$ be the quotient of this algebra,
 $\fH=H^\bu(G/B)\oplus H^\bu(G/B)/\k (\omega, -1)$,
where $\omega$, $1$  are the canonical generators of
$H^{\op{top}}(G/B)$, $H^0(G/B)$ respectively.

\begin{proposition}
The principal block of the center $\zl_0$ contains a canonically
defined subalgebra  isomorphic to $\fH$.
\end{proposition}

\proof The fibration $pr:\Nt\to G/B$ induces
 a short exact sequence of $G\times \C^*$ equivariant sheaves on
$\Nt$:
$$0\to T^{vert}_\Nt=z^{-2}\otimes pr^* T^*_{G/B} \to T_\Nt \to
T^{hor}_\Nt=pr^*T_{G/B} \to 0,$$ where $T^{vert}_\Nt$, $T^{hor}_\Nt$
are, respectively, the horizontal and the vertical tangent spaces.
In particular we get embeddings $$z^{-2i} \otimes pr^*\Omega^i_{G/B}\imbed
\Lambda^i T_{\Nt},$$
$$\Lambda^d(T^{vert}_\Nt)\otimes \Lambda
^i(pr^*(T_{G/B}))\imbed \Lambda ^{d+i}(T_\Nt),$$ where $d=\dim
(G/B)$. Notice that $\Lambda^d(T^{vert}_\Nt)\cong z^{-2d}\otimes
pr^*(\Omega^d_{G/B})$. Also a standard isomorphism
$\Lambda^d(V^*)\otimes \Lambda ^i(V)\cong \Lambda^{d-i}(V^*)$ for a
$d$-dimensional vector space $V$ yields $\Omega^d(G/B)\otimes
\Lambda^i(T_{G/B})\cong \Omega^{d-i}_{G/B}$. Thus the second
embedding above can be rewritten as: $$z^{-2d}\otimes
pr^*\Omega^{d-i}_{G/B}\imbed \Lambda ^{d+i}T_{\Nt}.$$

In fact, the two embeddings are easily seen to give isomorphisms of
sheaves on $G/B$: $$\Lambda^i(T_\Nt)^{-2i}\cong \Omega^i_{G/B}\cong
\Lambda^{2d-i}(T_\Nt)^{-2d}.$$

It is clear that the subsheaf
$\oplusl_i\Lambda^i(T_\Nt)^{-2i}\subset pr_*(\Lambda^\bu(T_\Nt))$ is
a  subsheaf of subalgebras isomorphic to $\Omega^\bu_{G/B}$. Thus, using
Corollary \ref{center} we get a subalgebra $\oplusl_i
H^i(\Lambda^i(T_\Nt)^{-2i})\subset \zl_0$, isomorphic to
$H^\bu (G/B)=\oplusl H^i(\Omega^i_{G/B})$.

 Similarly we get a
subsheaf $\oplusl_i \Lambda^{d+i}(T_\Nt)^{-2d}\subset
pr_*(\Lambda^\bu(T_\Nt))$ with zero multiplication; multiplication
of this sheaf by the sheaf $\Lambda^i(T_\Nt)^{-2i}$ also vanishes
for $i>0$. Thus we get a subspace $$H^\bu(\oplusl_i
\Lambda^{2d-i}(T_\Nt)^{-2d})\cong \oplusl H^i(\Omega^i_{G/B})=H^\bu(G/B)
\subset\zl_0.$$ Together with previously constructed subalgebra
$H^\bu (\oplusl_i\Lambda^i(T_\Nt)^{-2i})$ it clearly generates a
subalgebra canonically isomorphic to $\fH$. $\square$

\medskip

A $|W|$-dimensional subalgebra $\underline{\zl_0} \subset \zl_0$
isomorphic to $H^\bu(G/B)$ was described in \cite{BG}, where it was
obtained by ramification of the center of the De Concini-Kac quantum
algebra $\fU$. The subalgebra $\underline{\zl_0}$ can be defined as
the intersection of $Z_{HC} \otimes_{Z_l \cap Z_{HC}} \k$ with the
principal block $\ul_0$ and is usually referred to as the
Harish-Chandra part of the center of $\ul_0$.

In case when the root system of $\g$ is simple and simply laced, 
another $|W|$-dimensional subspace $\overline{\zl_0} \subset \zl_0$
was constructed in \cite{La} using the quantum Fourier transform 
as defined in 
\cite{LM}.  The subspace  $\overline{\zl_0}$ is shown to be 
an ideal in $\zl_0$ with zero multiplication, such that
$\op{Nilrad}(\zl_0)\cdot \overline{\zl_0} =0$. The 
 intersection $\underline{\zl_0} \cap
\overline{\zl_0}$ is one-dimensional and the subalgebra 
$\underline{\zl_0} +\overline{\zl_0} \subset \zl_0$ is isomorphic to $\fH$.

\bigskip {\bf Acknowledgements.} The first-named author would like
to thank Vladimir Drinfeld and Dennis Gaitsgory for useful
discussions. In particular, the idea of \S 2.4 stems from Drinfeld's ideas
on DG-schemes arising in the local geometric Langlands duality; it has been
independently discovered by Gaitsgory.

R.B. was partially supported by NSF, DARPA and the Sloan Foundation.
A.L. is grateful to the Max-Plank-Institut f\"ur Mathematik in Bonn
for financial support and hospitality.

\end{document}